\documentclass[amsmath,amssymb,prl,hyperlink,twocolumn]{revtex4}

	\usepackage{graphicx}
	\usepackage{soul}
	\usepackage[colorlinks=true,citecolor=blue,linkcolor=magenta]{hyperref}
	\usepackage[usenames]{color}
	\usepackage{amsfonts}
	\usepackage{color}
	\usepackage{booktabs}
	\usepackage{multirow}
	\usepackage{float}
    \usepackage{times}
    \usepackage[english]{babel}

\begin{document}

\preprint{APS/123-QED}

\title{Deep Representation Learning 
\\ for Dynamical Systems Modeling}

\author{Anna Shalova}
\email{Anna.Shalova@skoltech.ru}
\affiliation{Skolkovo Institute of Science and Technology
}%
\author{Ivan Oseledets}
\email{I.Oseledets@skoltech.ru}
\affiliation{Skolkovo Institute of Science and Technology
}%

\date{\today}

\begin{abstract}
Proper states' representations are the key to the successful dynamics modeling of chaotic systems. Inspired by recent advances of deep representations in various areas such as natural language processing and computer vision, we propose the adaptation of the state-of-art Transformer model in application to the dynamical systems modeling. The model demonstrates promising results in trajectories generation as well as in the general attractors' characteristics approximation, including states' distribution and Lyapunov exponent.

\end{abstract}

\maketitle

Revealing connections between the neural networks and dynamical systems is a long-standing quest for computational science. It is essential both for the usage of neural networks in application to the physical modeling and explanation of the deep neural models' behavior in general. Such approaches as neural ordinary differential equations \cite{chen2018odenet} and physical-informed neural networks \cite{raissi2019pinn} directly exploit connections between neural models and differential equations. Moreover, investigating neural networks' potential from the dynamical systems perspective has recently drawn a lot of attention \cite{battista2020capacity, chang2019antisymmetricrnn, laurent2016recurrentchaos}. In this work, we address the inverse problem and demonstrate the ability of the  probabilistic neural model \emph{transformer} \cite{vaswani2017transformer} to restore the chaotic attractors' structure. Given trajectories of a dynamical system of the form:
\begin{equation}\label{dynGPT-2:sys}
  \frac{dx}{dt} = f(x), \quad x(0)=x_0, \quad x \in \mathbb{R}^d,
\end{equation}
for different $x_0$ and sufficiently long time range, we would like to learn a model that would be able to predict the dynamics of such a system. 

Particular attention is accorded to the chaotic systems whose characteristic feature is the sensitivity to initial conditions \cite{ott2002chaos}. The exponential divergence of their trajectories leads to the instability of the approximation models. Classical approaches include dynamic mode decomposition \cite{kutz2016dmd, brunton2016dmd}, Koopman operator approach \cite{lusch2018deepkoopman, yeung2019deepkoopman}, reservoir computing \cite{pathak2018reservoir, jiang2019reservoir2}, and other methods based on the representation learning in the pre-determined space. With the rise of neural network models, a lot of attempts were made to model the dynamics of the systems in deep representations space using different neural architectures, including recurrent models (RNN) \cite{yu2017ttrnn} and residual networks (ResNet) \cite{chashchin2019resnet}. But despite being powerful interpolation models, the mentioned methods are not able to preserve trajectories from leaving the invariant set of an attractor. 

Chaotic systems are extremely challenging because trajectories are fundamentally unpredictable for the time exceeding the divergence time for the given numerical precision. The alternative approach is to approximate the probability distribution of the system's trajectories, which also allows the learning of informative deep states representations.

In this paper, we propose a new approach for black-box learning of dynamical systems from observations. Our inspiration comes from recent breakthrough advances in Natural Language Processing (NLP), where so-called \emph{transformers} \cite{vaswani2017transformer} have shown to be very effective for discrete sequence modeling. In contrast to NLP tasks, our state $x$ is continuous, so we first discretize our trajectory $x(t)$ using finite time step $\tau,$ and obtain time series $x_1, \ldots, x_T, \quad x_k = x(k \tau)$. Secondly, we discretize the trajectories in space: we embed the set $X$ of all possible values of $x_k$ into a $d$-dimensional cube and introduce a uniform  grid in each dimension with $n$ points. Then, each $x_k$ can be represented as an integer from $1$ to $n^d$, which is our discrete vocabulary. 

Once we have discrete trajectories, we propose to use the GPT-2 model \cite{radford2019gpt}. 
This model takes as an input the first $K$ states of the trajectory and reconstructs the remaining ones.
The GPT-2 model became a breakthrough in language modeling that allowed capturing of more complex and long-term dependencies than and other deep neural network models, such as RNNs or ResNet \cite{he2016resnet}. Transformers has also been successfully applied to the quantum states modeling \cite{carrasquilla2019quantum} and video representations learning \cite{sun2019contrastive}. We propose to use GPT-2 to learn deep states' representations of the dynamical systems as well as approximate the dynamics in the discrete space. We evaluate the approach from different perspectives, including the trajectories generation and approximation of the general properties of the systems.  

 We are given $M$ trajectories of length $T$, $x^{(i)}(t_k), i = 1, \ldots, M, k = 1, \ldots, T$ which arise from the dynamical system \eqref{dynGPT-2:sys} by sampling different initial conditions $x_0$, and then taking $x(t)$ at discrete times  $t_k = k \tau$. Our goal is to train an autoregressive model that given the trajectory $x_1, \ldots, x_K$ predicts its next state $x_{K+1} = x(t_{k+1})$. 

We suggest modeling the evolution of the dynamical system in the discrete latent space instead of the continuous one. By reducing the spatial resolution, we hope to get better long-term behavior. As an \emph{encoder} we use the most straightforward possible mapping: equidistant grid discretization with $n$ segments in every dimension, which leads to the $n^d$ size of the ``vocabulary'' to represent the state. As a \emph{decoder} we can use piecewise constant continuation, i.e.,  to decode discrete states into the continuous space, we use the center of the cube corresponding to the current class.

\begin{figure}
    \centering
    \includegraphics[width= 0.5\textwidth]{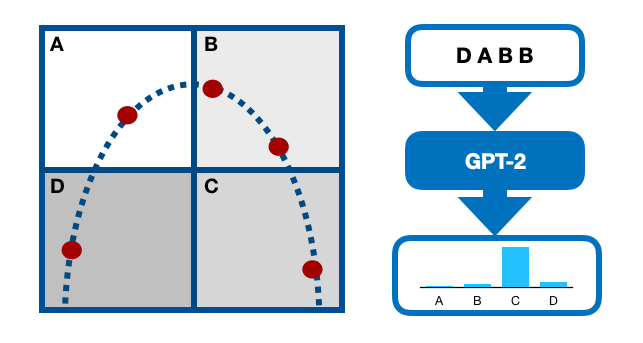}
    \caption{Scheme of the proposed system}
    \label{fig:scheme}
\end{figure}

An example of the proposed system is presented in Figure \ref{fig:scheme}. Transition to the discrete space allows us to reduce the regression task to the classification problem with the number of classes equal to the size of the vocabulary. Given the trajectory $\{i_0, i_1, \ldots, i_T\}$, the model is trained to maximize the log-probability of the next state conditioned by the previous ones $p(i_t | i_0, \ldots, i_{T-1})$ for the observed trajectory. As a result, it learns conditional factorization of the joint probability of the given sequences. This is exactly the task that is known in NLP as language modeling, and we propose to solve it using the GPT-2 model that is based on the self-attention mechanism.

Self-attention was initially presented in the \emph{transformer} model in \cite{vaswani2017transformer}. The transformer is a sequence to sequence model that, in the NLP task setting,  learns a mapping from one text to another. Such models traditionally consist of encoder and decoder, where encoder learns some representation of an initially given sequence, and decoder learns dynamics in the target space conditioned by the given representation. Lately, it was proposed to use encoder and decoder separately, and the GPT-2 is exactly an unconditioned version of the transformer decoder \cite{radford2019gpt}.

In NLP text is usually viewed as a sequence of the word- or character-level language units that are also called \emph{tokens}. Most of the modern models in natural language learn dependencies in the space of vector representations of the tokens. 
An intuitive approach is to measure relations between tokens representations as a scalar product of their embeddings. Similar to texts, trajectories of dynamical systems can be viewed as a sequence of discrete states from some vocabulary. 

Encoder and decoder of the transformer have a similar structure with the only difference that the decoder has auto-regressive nature. On each step, the decoder operates only with the previous tokens while the encoder sees the whole input sequence. Although it is not important for representation learning, holding causality is crucial for dynamics modeling. 

GPT-2 model is a sequence of self-attention layers with a multiclass-classification head with the number of classes equal to the size of the vocabulary. The self-attention layer takes a sequence of tokens representations as an input and provides an output of the same shape. The layer has several parallel attention-heads; their outputs are concatenated and projected to get the output of the layer. Each head has three different projections, two of them are used to calculate the closeness of tokens, and then the third representations are summed with the weights determined by the obtained similarity. 
A more detailed explanation of the self-attention mechanism can be found in the original paper \cite{vaswani2017transformer}. 

Another feature of transformer models is the usage of the positional embeddings. They are constant vectors that are added to the token embeddings, which helps self-attention layers to determine the position of each token in the given sequence.

The GPT-2 model is trained to minimize the cross-entropy loss that is standard for classification models. It corresponds to the maximization of conditional log-probabilities of every token from the given trajectories:
$$\mathcal{L} = -\frac1{N}\sum_i\log p(x_i|x_1, .. x_{i - 1}).$$

It is equivalent to the maximization of KL-divergence between the real conditional distribution of the training set and modeled conditional distributions.
The input of the model is limited, which means that we can generate trajectory conditioned by up to 1023 previous states. Similar to the text generation, trajectories are predicted sequentially $i_t = f_{t-1}( i_0, \ldots, i_{t-1})$. On each step the model predicts discrete probability distribution containing probabilities of the next state to belong to each cube:

$$p(i_t) = \text{softmax}(z_t/T),$$
where T (also called temperature) is a normalizing constant used before the softmax layer. Temperature is a trade-off between the variability and quality of the generated sequences; it allows balancing between the probability of the predicted states and the coverage of the trajectories distribution. In our experiments, we used $T \to 0$, which corresponds to the choice of the most probable token on each step. We also tested higher temperature sampling for the distribution approximation, but have seen no significant advantage.

In all our experiments, we use the GPT-2 model consisting of 12 self-attention layers and representations dimension 768. The input of the model is limited by 1024 tokens, which means that it looks up to 1024 steps back that allows capturing long-term dependencies. We use PyTorch implementation of GPT-2 as well as AdamW optimizer from \cite{Wolf2019HuggingFacesTS}. 
We use default parameters of the optimizer except for the initial learning rate = 0.00005.

We demonstrate GPT-2 performance from different perspectives for the following dynamical systems:
\begin{enumerate}
\item \textbf{Lorenz attractor}
is a set of chaotic solutions of the Lorenz system that describes two-dimensional fluids flow:
\begin{equation}\label{dynsys:lorentz}
  \begin{split}
    \frac{dx}{dt} = \sigma(y - x), 
    \frac{dy}{dt} = x(\rho - z) - y, 
    \frac{dz}{dt} = xy - \beta z.
  \end{split}
\end{equation}
In our experiment we use parameters $\sigma = 28$, $\rho = 10 $, $\beta = 2.67$. Our training and test set consist of 1000 and 20 trajectories with initial states generated from uniform distribution [-0.1, 0.1]. For each trajectory we run simulation for $t = 1000$ with a time step $\tau = 0.03$.
\item \textbf{Rossler attractor} is an attractor for the Rossler system:
\begin{equation}\label{dynsys:rossler}
  \begin{split}
    \frac{dx}{dt} = -y - z, \frac{dy}{dt} = x+ay, \frac{dz}{dt} = b + z(x - c).
  \end{split}
\end{equation}
We use Rossler attractor with $a = 0.15, b =0.2 , c = 10 $. Similar to the Lorenz attractor we have 1000 and 20 trajectories in train and test sets. Each trajectory corresponds to the simulation for $t = 1000$ and $\tau = 0.1$. Initial states are sampled from distribution $x = x_0 + \epsilon$ with $x_0 = [5, 0, 0]$ and uniform $\epsilon \sim Uni(-1, 1)$.
\item \textbf{Henon map} is a discrete-time dynamical system:

\begin{equation}\label{dynsys:henon}
  \begin{split}
    x_{n + 1} =1 - ax_n^2 + y_n, y_{n + 1} = bx_n. \\
  \end{split}
\end{equation}
We use map with $a = 1.4, b = 0.3$. As the system is two-dimensional, we have only 100 and 10 trajectories in train and test sets. Each trajectory corresponds to the simulation for $t = 10000$ steps. Initial states are sampled from distribution $x = x_0 + \epsilon$ with $x_0 = [-0.95, 0.35]$ and uniform $\epsilon \sim Uni(-0.05, 0.05)$.
\end{enumerate}

The distinctive feature of the proposed model is that predictions are limited by the set of the states presented in the training data. 
This fact does not guarantee that predicted trajectories will have desirable structure but allows simple evaluation of the distribution approximation.

There are a lot of tasks for which modeling of the general behavior of the attractor is not less important that precise trajectories reconstruction; for example, estimation of the Lyapunov exponents. 

We compare empirical states distributions obtained from the trajectories in the discrete space. We use $10^5$-$10^6$ states in each case, which is shown to be sufficient for the Lyapunov exponent approximation. As a metric between distributions, we use the Wasserstein distance. It can be interpreted a minimal work required to transform one distribution to another. For the discrete distributions $u$ and $v$ defined in the same metric space Wasserstein distance $W(u, v)$ is a solution of an optimal transport problem.

Results are presented in Table \ref{tab:wd}. For the Henon map and Rossler attractor we use only the grid with the smallest discretization step (N = 50). For the Lorenz system, we investigate the influence of the grid size on the Wasserstein distance. Note that the metric itself depends on the discretization, so we provide a reference value - the Wasserstein distance between two batches of the simulated trajectories. In all cases, Wasserstein distance between true and generated distributions is of the order as the distance between two true distributions. This means that regardless of the discretization size, the model is approximately equally effective. Different grid sizes also correspond to the different ratio of time and space characteristics of the system, which indicates how many cubes the model is expected to pass in one step. 

\begin{table}[h]
\caption{ \label{tab:wd} Wasserstein distance of the generated states' distributions}
\begin{tabular}{|l|l|l|l|l|l|}
\hline
             & \multicolumn{3}{c|}{\textbf{Lorenz}} & \textbf{Henon} & \textbf{Rossler} 
\\ \hline
            & N = 20 &  N = 35 & N = 50 & N = 50 & N = 50
\\\hline
GPT-2 vs True1& 3.0379 & 2.4122   & 2.1777  & 0.00499 & 0.0557
\\\hline
True1 vs True2 & 3.1203 & 3.1134  & 2.0274  & 0.00379 & 0.0822
\\\hline
\end{tabular}
\end{table}

\label{sec:exponent}
From the obtained states' distribution, we can also evaluate the largest Lyapunov exponent. Lyapunov spectrum in general and the largest Lyapunov exponent in particular are one of the distinct characteristics of the chaotic system. It can be called a measure of the chaos of the system. If initial states of two trajectories in the phase space differ by $\delta x_0 \to 0$, their divergence is expected to be 
$$\delta x (t) \sim e^{\lambda t}\delta x_0,$$
where $\lambda$ is the Lyapunov exponent. Depending on the orientation of the difference, there could be several values of $\lambda$ called the spectrum of Lyapunov exponents. The largest Lyapunov exponent can be defined as:
$$\lambda = \lim_{t\to\infty}\lim_{\delta x_0 \to 0}\frac1{t}\ln{\frac{\delta x(t)}{\delta x_0}}.$$

In this experiment, we use one of the approaches to the maximum Lyapunov exponent estimation that, in our case, relies on the states' distribution \cite{aston1999exponent}. For the discrete-time dynamical system $x_{n + 1} = F(x_n)$ the largest Lyapunov exponent can be calculated using time averaging:

$$\lambda = \lim_{N \to \infty} \frac1{N} \log \|\prod_{i = 1}^NF'(x_i)\|  = \lim_{N \to \infty} \lambda_n,$$
where $\|F'(x_i)\|$ is a spectral norm of the Jacobian in state $x_i$. In \cite{aston1999exponent} it was proposed to model $\lambda_n(n)$ in the following form:
$$\lambda_n(n) = \lambda + \frac{c_1}{n} + \frac{c_2}{n^2} + O(\frac{c_2}{n^3})$$
and estimate $\lambda$ from a small number of first elements of the sequence $\{\lambda_n\}$. 

We approximate the largest Lyapunov exponent of the Henon map. In this experiment, we consider Jacobian for every state known and equal to the Jacobian of the system calculated in the center of the corresponding cube. We calculate $\{\lambda_i: 1\leq i \leq15\}$ from ten trajectories of 10000 states generated by GPT-2 conditioned by 100 states each.

As a result we got $\lambda_{GPT-2} =  0.4161$ while known estimated value is $\lambda = 0.4192$ \cite{aston1999exponent}. We can see that the Lyapunov exponent approximated from the generated trajectories is close to its real value, which demonstrates that we have a sufficiently good approximation of conditional distributions with the lengths of conditioning up to 14.

We also propose alternative approaches to the trajectories divergence estimation. The chaotic structure of the initial model does not allow a precise evaluation of the generated trajectories.  
During discretization, we lose micro-scale information about trajectories; as a result, one trajectory in the discrete space can correspond not to the one but the range of trajectories in the continuous space. So the comparison of the restored trajectory to the unique continuation is not correct.

As a solution to this problem, we suggest the following algorithm. Having one trajectory in the continuous space with initial state $x_0$ and its mapping to the discrete space, we want to find as many possible continuations as possible and estimate the divergence time using the closest one to the sequence generated by our model.

\begin{figure}
\includegraphics[width = 0.5\textwidth]{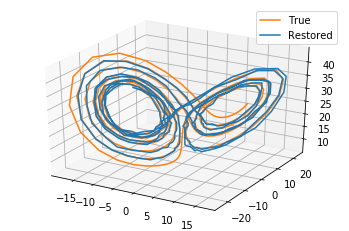}
\caption{\label{fig:traj} Simulated and Restored by GPT-2 trajectories of Lorenz attractor.}
\end{figure}

If the trajectory is conditioned by $k$ states, from the definition of the largest Lyapunov exponent the difference between the initial states can be estimated as:
$$\delta x_0 \sim l_d e^{-\lambda k \tau}  $$
where $l_d$ is the radius of the sphere wrapped around the cell of the discretization grid. So it can be expected that trajectories of $k$ states with initial condition from the normal distribution $\mathcal{N}(x_0, \delta x_0 I)$ will have the same mapping to the discrete space.

For every sampled trajectory, we check that the first $k$ representations match the given condition and then find the time for which it matches the predicted trajectory. As a result, we use the maximum time among all sampled trajectories. Still, due to the chaotic nature of the system, we have no guarantees that there are no longer generated trajectories corresponding to some true trajectories. 

We conduct experiments for both Lorenz and Rossler attractors. We use $k = 100$ to limit the area of initial states and provide a sufficient acceptance rate of the sampled trajectories. We also use models with the smallest grid $N = 50$. Examples of the trajectories generated by GPT-2 for Lorenz systems are shown in Figure \ref{fig:traj}. It can be noticed that even after divergence from the reference trajectories, continuations generated by the GPT-2 have a similar structure and are not distinguishable from true samples.

The characteristic time of the system can be evaluated as inverse maximum Lyapunov exponent $t_d = 1/\lambda$. As it was shown in our experiment, time-discrete transformer models have very close time scales to the original systems. The estimated divergence time is 1.11 for the Lorenz attractor and 4.59 for the Rossler attractor, while inverse Lyapunov exponents are 1.04 and 5.33 respectfully \cite{wolf1985lyapTS}. It demonstrates that modeled time-discrete continuations diverge from the simulated trajectories on the same time scale as simulated trajectories diverge from each other.  

\begin{figure}
    \centering
    \includegraphics[width = 0.45\textwidth]{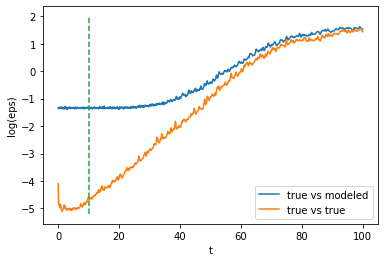}
    \caption{Expected difference between trajectories from time. Dotted line indicates end of conditioning for the GPT-2 model.}
    \label{fig:diff}
\end{figure}

For the Rossler attractor, we also compare the dependence of the expected difference between trajectories from time. Firstly we randomly select 300 trajectories of length 1000 from the test dataset. Then, for every sample, we generate a trajectory with GPT-2 using its first 100 states as a condition. We also simulate a paired trajectory for every sample so that discretizations of their first 100 steps are the same. That allows us to estimate divergence between two real trajectories. The result is presented in Figure \ref{fig:diff}. We can see that the long-term behavior of two curves is close, including the slope corresponding to the Lyapunov exponent. That additionally proves that the GPT-2 is able to reproduce the attractor behavior.

To sum up, we demonstrated that the GPT-2 model could be successfully adopted for dynamical systems modeling. With a simple discretization technique, it is able to restore the general structure of attractors in terms of states' distributions, generated trajectories, and divergence characteristics. Despite the limited accuracy of single prediction, our model is more efficient for a long-term trajectories generation, and, in contrast to the regression models, the probabilistic approach guarantees the existence of the invariant state set.
\bigskip

\emph{Acknowledgements.} This work was supported by Ministry of Higher Education and Science of the Russian Federation under grant 14.756.31.0001.

\bibliography{apssamp.bib}

\begin{thebibliography}{22}
\expandafter\ifx\csname natexlab\endcsname\relax\def\natexlab#1{#1}\fi
\expandafter\ifx\csname bibnamefont\endcsname\relax
  \def\bibnamefont#1{#1}\fi
\expandafter\ifx\csname bibfnamefont\endcsname\relax
  \def\bibfnamefont#1{#1}\fi
\expandafter\ifx\csname citenamefont\endcsname\relax
  \def\citenamefont#1{#1}\fi
\expandafter\ifx\csname url\endcsname\relax
  \def\url#1{\texttt{#1}}\fi
\expandafter\ifx\csname urlprefix\endcsname\relax\def\urlprefix{URL }\fi
\providecommand{\bibinfo}[2]{#2}
\providecommand{\eprint}[2][]{\url{#2}}

\bibitem[{\citenamefont{Chen et~al.}(2018)\citenamefont{Chen, Rubanova,
  Bettencourt, and Duvenaud}}]{chen2018odenet}
\bibinfo{author}{\bibfnamefont{T.~Q.} \bibnamefont{Chen}},
  \bibinfo{author}{\bibfnamefont{Y.}~\bibnamefont{Rubanova}},
  \bibinfo{author}{\bibfnamefont{J.}~\bibnamefont{Bettencourt}},
  \bibnamefont{and} \bibinfo{author}{\bibfnamefont{D.~K.}
  \bibnamefont{Duvenaud}}, in \emph{\bibinfo{booktitle}{Advances in neural
  information processing systems}} (\bibinfo{year}{2018}), pp.
  \bibinfo{pages}{6571--6583}.

\bibitem[{\citenamefont{Raissi et~al.}(2019)\citenamefont{Raissi, Perdikaris,
  and Karniadakis}}]{raissi2019pinn}
\bibinfo{author}{\bibfnamefont{M.}~\bibnamefont{Raissi}},
  \bibinfo{author}{\bibfnamefont{P.}~\bibnamefont{Perdikaris}},
  \bibnamefont{and} \bibinfo{author}{\bibfnamefont{G.~E.}
  \bibnamefont{Karniadakis}}, \bibinfo{journal}{Journal of Computational
  Physics} \textbf{\bibinfo{volume}{378}}, \bibinfo{pages}{686}
  (\bibinfo{year}{2019}).

\bibitem[{\citenamefont{Battista and Monasson}(2020)}]{battista2020capacity}
\bibinfo{author}{\bibfnamefont{A.}~\bibnamefont{Battista}} \bibnamefont{and}
  \bibinfo{author}{\bibfnamefont{R.}~\bibnamefont{Monasson}},
  \bibinfo{journal}{Physical Review Letters} \textbf{\bibinfo{volume}{124}},
  \bibinfo{pages}{048302} (\bibinfo{year}{2020}).

\bibitem[{\citenamefont{Chang et~al.}(2019)\citenamefont{Chang, Chen, Haber,
  and Chi}}]{chang2019antisymmetricrnn}
\bibinfo{author}{\bibfnamefont{B.}~\bibnamefont{Chang}},
  \bibinfo{author}{\bibfnamefont{M.}~\bibnamefont{Chen}},
  \bibinfo{author}{\bibfnamefont{E.}~\bibnamefont{Haber}}, \bibnamefont{and}
  \bibinfo{author}{\bibfnamefont{E.~H.} \bibnamefont{Chi}},
  \bibinfo{journal}{arXiv preprint arXiv:1902.09689}  (\bibinfo{year}{2019}).

\bibitem[{\citenamefont{Laurent and von
  Brecht}(2016)}]{laurent2016recurrentchaos}
\bibinfo{author}{\bibfnamefont{T.}~\bibnamefont{Laurent}} \bibnamefont{and}
  \bibinfo{author}{\bibfnamefont{J.}~\bibnamefont{von Brecht}},
  \bibinfo{journal}{arXiv preprint arXiv:1612.06212}  (\bibinfo{year}{2016}).

\bibitem[{\citenamefont{Vaswani et~al.}(2017)\citenamefont{Vaswani, Shazeer,
  Parmar, Uszkoreit, Jones, Gomez, Kaiser, and
  Polosukhin}}]{vaswani2017transformer}
\bibinfo{author}{\bibfnamefont{A.}~\bibnamefont{Vaswani}},
  \bibinfo{author}{\bibfnamefont{N.}~\bibnamefont{Shazeer}},
  \bibinfo{author}{\bibfnamefont{N.}~\bibnamefont{Parmar}},
  \bibinfo{author}{\bibfnamefont{J.}~\bibnamefont{Uszkoreit}},
  \bibinfo{author}{\bibfnamefont{L.}~\bibnamefont{Jones}},
  \bibinfo{author}{\bibfnamefont{A.~N.} \bibnamefont{Gomez}},
  \bibinfo{author}{\bibfnamefont{{\L}.}~\bibnamefont{Kaiser}},
  \bibnamefont{and}
  \bibinfo{author}{\bibfnamefont{I.}~\bibnamefont{Polosukhin}}, in
  \emph{\bibinfo{booktitle}{Advances in neural information processing systems}}
  (\bibinfo{year}{2017}), pp. \bibinfo{pages}{5998--6008}.

\bibitem[{\citenamefont{Ott}(2002)}]{ott2002chaos}
\bibinfo{author}{\bibfnamefont{E.}~\bibnamefont{Ott}},
  \emph{\bibinfo{title}{Chaos in dynamical systems}}
  (\bibinfo{publisher}{Cambridge university press}, \bibinfo{year}{2002}).

\bibitem[{\citenamefont{Kutz et~al.}(2016)\citenamefont{Kutz, Brunton, Brunton,
  and Proctor}}]{kutz2016dmd}
\bibinfo{author}{\bibfnamefont{J.~N.} \bibnamefont{Kutz}},
  \bibinfo{author}{\bibfnamefont{S.~L.} \bibnamefont{Brunton}},
  \bibinfo{author}{\bibfnamefont{B.~W.} \bibnamefont{Brunton}},
  \bibnamefont{and} \bibinfo{author}{\bibfnamefont{J.~L.}
  \bibnamefont{Proctor}}, \emph{\bibinfo{title}{Dynamic mode decomposition:
  data-driven modeling of complex systems}} (\bibinfo{publisher}{SIAM},
  \bibinfo{year}{2016}).

\bibitem[{\citenamefont{Brunton et~al.}(2016)\citenamefont{Brunton, Proctor,
  and Kutz}}]{brunton2016dmd}
\bibinfo{author}{\bibfnamefont{S.~L.} \bibnamefont{Brunton}},
  \bibinfo{author}{\bibfnamefont{J.~L.} \bibnamefont{Proctor}},
  \bibnamefont{and} \bibinfo{author}{\bibfnamefont{J.~N.} \bibnamefont{Kutz}},
  \bibinfo{journal}{Proceedings of the national academy of sciences}
  \textbf{\bibinfo{volume}{113}}, \bibinfo{pages}{3932} (\bibinfo{year}{2016}).

\bibitem[{\citenamefont{Lusch et~al.}(2018)\citenamefont{Lusch, Kutz, and
  Brunton}}]{lusch2018deepkoopman}
\bibinfo{author}{\bibfnamefont{B.}~\bibnamefont{Lusch}},
  \bibinfo{author}{\bibfnamefont{J.~N.} \bibnamefont{Kutz}}, \bibnamefont{and}
  \bibinfo{author}{\bibfnamefont{S.~L.} \bibnamefont{Brunton}},
  \bibinfo{journal}{Nature communications} \textbf{\bibinfo{volume}{9}},
  \bibinfo{pages}{1} (\bibinfo{year}{2018}).

\bibitem[{\citenamefont{Yeung et~al.}(2019)\citenamefont{Yeung, Kundu, and
  Hodas}}]{yeung2019deepkoopman}
\bibinfo{author}{\bibfnamefont{E.}~\bibnamefont{Yeung}},
  \bibinfo{author}{\bibfnamefont{S.}~\bibnamefont{Kundu}}, \bibnamefont{and}
  \bibinfo{author}{\bibfnamefont{N.}~\bibnamefont{Hodas}}, in
  \emph{\bibinfo{booktitle}{2019 American Control Conference (ACC)}}
  (\bibinfo{organization}{IEEE}, \bibinfo{year}{2019}), pp.
  \bibinfo{pages}{4832--4839}.

\bibitem[{\citenamefont{Pathak et~al.}(2018)\citenamefont{Pathak, Hunt, Girvan,
  Lu, and Ott}}]{pathak2018reservoir}
\bibinfo{author}{\bibfnamefont{J.}~\bibnamefont{Pathak}},
  \bibinfo{author}{\bibfnamefont{B.}~\bibnamefont{Hunt}},
  \bibinfo{author}{\bibfnamefont{M.}~\bibnamefont{Girvan}},
  \bibinfo{author}{\bibfnamefont{Z.}~\bibnamefont{Lu}}, \bibnamefont{and}
  \bibinfo{author}{\bibfnamefont{E.}~\bibnamefont{Ott}},
  \bibinfo{journal}{Physical review letters} \textbf{\bibinfo{volume}{120}},
  \bibinfo{pages}{024102} (\bibinfo{year}{2018}).

\bibitem[{\citenamefont{Jiang and Lai}(2019)}]{jiang2019reservoir2}
\bibinfo{author}{\bibfnamefont{J.}~\bibnamefont{Jiang}} \bibnamefont{and}
  \bibinfo{author}{\bibfnamefont{Y.-C.} \bibnamefont{Lai}},
  \bibinfo{journal}{Physical Review Research} \textbf{\bibinfo{volume}{1}},
  \bibinfo{pages}{033056} (\bibinfo{year}{2019}).

\bibitem[{\citenamefont{Yu et~al.}(2017)\citenamefont{Yu, Zheng, Anandkumar,
  and Yue}}]{yu2017ttrnn}
\bibinfo{author}{\bibfnamefont{R.}~\bibnamefont{Yu}},
  \bibinfo{author}{\bibfnamefont{S.}~\bibnamefont{Zheng}},
  \bibinfo{author}{\bibfnamefont{A.}~\bibnamefont{Anandkumar}},
  \bibnamefont{and} \bibinfo{author}{\bibfnamefont{Y.}~\bibnamefont{Yue}},
  \bibinfo{journal}{arXiv preprint arXiv:1711.00073}  (\bibinfo{year}{2017}).

\bibitem[{\citenamefont{Chashchin et~al.}(2019)\citenamefont{Chashchin,
  Botchev, Oseledets, and Ovchinnikov}}]{chashchin2019resnet}
\bibinfo{author}{\bibfnamefont{A.}~\bibnamefont{Chashchin}},
  \bibinfo{author}{\bibfnamefont{M.}~\bibnamefont{Botchev}},
  \bibinfo{author}{\bibfnamefont{I.}~\bibnamefont{Oseledets}},
  \bibnamefont{and}
  \bibinfo{author}{\bibfnamefont{G.}~\bibnamefont{Ovchinnikov}},
  \bibinfo{journal}{arXiv preprint arXiv:1910.05233}  (\bibinfo{year}{2019}).

\bibitem[{\citenamefont{Radford et~al.}(2019)\citenamefont{Radford, Wu, Child,
  Luan, Amodei, and Sutskever}}]{radford2019gpt}
\bibinfo{author}{\bibfnamefont{A.}~\bibnamefont{Radford}},
  \bibinfo{author}{\bibfnamefont{J.}~\bibnamefont{Wu}},
  \bibinfo{author}{\bibfnamefont{R.}~\bibnamefont{Child}},
  \bibinfo{author}{\bibfnamefont{D.}~\bibnamefont{Luan}},
  \bibinfo{author}{\bibfnamefont{D.}~\bibnamefont{Amodei}}, \bibnamefont{and}
  \bibinfo{author}{\bibfnamefont{I.}~\bibnamefont{Sutskever}},
  \bibinfo{journal}{OpenAI Blog} \textbf{\bibinfo{volume}{1}},
  \bibinfo{pages}{9} (\bibinfo{year}{2019}).

\bibitem[{\citenamefont{He et~al.}(2016)\citenamefont{He, Zhang, Ren, and
  Sun}}]{he2016resnet}
\bibinfo{author}{\bibfnamefont{K.}~\bibnamefont{He}},
  \bibinfo{author}{\bibfnamefont{X.}~\bibnamefont{Zhang}},
  \bibinfo{author}{\bibfnamefont{S.}~\bibnamefont{Ren}}, \bibnamefont{and}
  \bibinfo{author}{\bibfnamefont{J.}~\bibnamefont{Sun}}, in
  \emph{\bibinfo{booktitle}{Proceedings of the IEEE conference on computer
  vision and pattern recognition}} (\bibinfo{year}{2016}), pp.
  \bibinfo{pages}{770--778}.

\bibitem[{\citenamefont{Carrasquilla et~al.}(2019)\citenamefont{Carrasquilla,
  Luo, P{\'e}rez, Milsted, Clark, Volkovs, and
  Aolita}}]{carrasquilla2019quantum}
\bibinfo{author}{\bibfnamefont{J.}~\bibnamefont{Carrasquilla}},
  \bibinfo{author}{\bibfnamefont{D.}~\bibnamefont{Luo}},
  \bibinfo{author}{\bibfnamefont{F.}~\bibnamefont{P{\'e}rez}},
  \bibinfo{author}{\bibfnamefont{A.}~\bibnamefont{Milsted}},
  \bibinfo{author}{\bibfnamefont{B.~K.} \bibnamefont{Clark}},
  \bibinfo{author}{\bibfnamefont{M.}~\bibnamefont{Volkovs}}, \bibnamefont{and}
  \bibinfo{author}{\bibfnamefont{L.}~\bibnamefont{Aolita}},
  \bibinfo{journal}{arXiv preprint arXiv:1912.11052}  (\bibinfo{year}{2019}).

\bibitem[{\citenamefont{Sun et~al.}(2019)\citenamefont{Sun, Baradel, Murphy,
  and Schmid}}]{sun2019contrastive}
\bibinfo{author}{\bibfnamefont{C.}~\bibnamefont{Sun}},
  \bibinfo{author}{\bibfnamefont{F.}~\bibnamefont{Baradel}},
  \bibinfo{author}{\bibfnamefont{K.}~\bibnamefont{Murphy}}, \bibnamefont{and}
  \bibinfo{author}{\bibfnamefont{C.}~\bibnamefont{Schmid}},
  \bibinfo{journal}{arXiv preprint arXiv:1906.05743}  (\bibinfo{year}{2019}).

\bibitem[{\citenamefont{Wolf et~al.}(2019)\citenamefont{Wolf, Debut, Sanh,
  Chaumond, Delangue, Moi, Cistac, Rault, Louf, Funtowicz
  et~al.}}]{Wolf2019HuggingFacesTS}
\bibinfo{author}{\bibfnamefont{T.}~\bibnamefont{Wolf}},
  \bibinfo{author}{\bibfnamefont{L.}~\bibnamefont{Debut}},
  \bibinfo{author}{\bibfnamefont{V.}~\bibnamefont{Sanh}},
  \bibinfo{author}{\bibfnamefont{J.}~\bibnamefont{Chaumond}},
  \bibinfo{author}{\bibfnamefont{C.}~\bibnamefont{Delangue}},
  \bibinfo{author}{\bibfnamefont{A.}~\bibnamefont{Moi}},
  \bibinfo{author}{\bibfnamefont{P.}~\bibnamefont{Cistac}},
  \bibinfo{author}{\bibfnamefont{T.}~\bibnamefont{Rault}},
  \bibinfo{author}{\bibfnamefont{R.}~\bibnamefont{Louf}},
  \bibinfo{author}{\bibfnamefont{M.}~\bibnamefont{Funtowicz}},
  \bibnamefont{et~al.}, \bibinfo{journal}{ArXiv}
  \textbf{\bibinfo{volume}{abs/1910.03771}} (\bibinfo{year}{2019}).

\bibitem[{\citenamefont{Aston and Dellnitz}(1999)}]{aston1999exponent}
\bibinfo{author}{\bibfnamefont{P.~J.} \bibnamefont{Aston}} \bibnamefont{and}
  \bibinfo{author}{\bibfnamefont{M.}~\bibnamefont{Dellnitz}},
  \bibinfo{journal}{Computer methods in applied mechanics and engineering}
  \textbf{\bibinfo{volume}{170}}, \bibinfo{pages}{223} (\bibinfo{year}{1999}).

\bibitem[{\citenamefont{Wolf et~al.}(1985)\citenamefont{Wolf, Swift, Swinney,
  and Vastano}}]{wolf1985lyapTS}
\bibinfo{author}{\bibfnamefont{A.}~\bibnamefont{Wolf}},
  \bibinfo{author}{\bibfnamefont{J.~B.} \bibnamefont{Swift}},
  \bibinfo{author}{\bibfnamefont{H.~L.} \bibnamefont{Swinney}},
  \bibnamefont{and} \bibinfo{author}{\bibfnamefont{J.~A.}
  \bibnamefont{Vastano}}, \bibinfo{journal}{Physica D: Nonlinear Phenomena}
  \textbf{\bibinfo{volume}{16}}, \bibinfo{pages}{285} (\bibinfo{year}{1985}).

\end{thebibliography}
\end{document}